\documentclass{amsart}
\usepackage{amsthm,amsfonts,amsmath,amscd,amssymb,latexsym,epsfig,rb}
\newcommand{\gl}{{\mathfrak g \mathfrak l}}

\newcommand{\su}{{\mathfrak s  \mathfrak u}}

\newcommand{\g}{{\mathfrak g}}         
\newcommand{\kk}{{\mathfrak k}}         
\newcommand{\cx}{{\mathbb C}}
\newcommand{\quat}{{\mathbb H}}
\newcommand{\proj}{{\mathbb P}}

\newcommand{\ad}{\operatorname{ad}}

\newcommand{\Lie}{\operatorname{Lie}}

\newcommand{\G}{{\mathcal G}}
\newcommand{\Gn}{{\mathcal G}_0}
\renewcommand{\A}{{\mathcal A}}
\newcommand{\AC}{{\mathcal A}^\cx}
\newcommand{\AH}{{\mathcal A}^\quat}
\newcommand{\ttt}{\tau_1,\tau_2,\tau_3}
\newcommand{\ttts}{{\tau_1,\tau_2,\tau_3;\sigma}}
\newcommand{\gk}{(\g,\kk)}

\newcommand{\thm}[2]{\textsc{#1:} #2}
\newcommand{\exercise}[1]{\thm{exercise}{#1}}
\newcommand{\claim}[1]{\thm{claim}{#1}}

\begin{document}

\title[Lie groups, Nahm's equations
 and hyperk\"ahler manifolds]{Lie groups, Nahm's equations\\
 and hyperk\"ahler manifolds}
\author[Roger Bielawski]{Lectures of 
Roger Bielawski\\
given at the summer school on Algebraic Groups}\thanks{Notes by Sven-S. Porst}
\date{G\"ottingen, July 2005}
\maketitle

The Nahm equations are a powerful tool for constructing hyperk\"ahler metrics on various algebraic manifolds, e.g. spaces of rational maps, coadjoint orbits, resolutions of Kleinian singularities.  In these lectures I shall concentrate on their relation to semisimple algebraic groups and adjoint orbits. The aim is to show that, on the one hand, the Nahm equations provide a powerful tool for constructing geometric structures on these objects and for explaining certain representation-theoretic puzzles (the Kostant-Sekiguchi correspondence \cite{Vergne}, action of the Weyl group on flag manifolds \cite{ABi}), and, on the other hand, to stress just how mysterious these geometric structures remain. This is so particularly for the case of hyperk\"ahler metrics on coadjoint orbits.
   \par
   There is very little in these lecture that is original: perhaps some extensions of known results (e.g. Theorem 3.4 or Corollary 3.5 in literature) and some of the approach. On the other hand, I allowed myself to speculate at some points.
   \par
   It is a pleasure to thank Victor Pidstrygach and Yuri Tschinkel for the invitation to give these lectures and Sven Porst for taking and typing up the notes.

\section{Quaternions and Lie groups}

We begin with an unusually complicated construction of a bi-invariant metric on a compact Lie group $G$ using the space of paths in its Lie algebra $\g$. This approach is extended to the K\"ahler and  hyperk\"ahler reductions in the complexified and quaternionised cases.

\subsection{Construct a bi-invariant metric}
Consider the sets of paths in a compact connected Lie group $G$ and its Lie algebra $\g$
\[
\A = \{ T:[0,1] \longrightarrow \g\}
\qquad \mathrm{and} \qquad
\G = \{g:[0,1] \longrightarrow G\}.
\]
$\A$ is an infinite-dimensional vector space and $\G$ is an infinite-dimensional Lie group which can be modelled by a Hilbert space. Denote by $\Gn$ be  the subgroup of loops starting at the unit element of $G$:
\[
\Gn = \{g\in\G ; g(0) = 1 = g(1)\} \subset \G
\]
Elements $g\in\G$ act on $\A$ by $(g.T)(s) = g(s)T(s)g(s)^{-1} - \dot{g}(s)g^{-1}(s)$ where  $\dot{g}$ is the derivative of $g$ at $0$. Thus $\A$ is the space of connections on the trivial $G$-bundle over $[0,1]$ and $\G$ is the gauge group. We shall drop the parameter $s$ from notation whenever possible, stating the action as:
\[
g.T = gTg^{-1} - \dot{g}g^{-1}.
\]
The action of  $\Gn$ on $\A$ is free. For any $T\in\A$ we have the set $\{g\in \G ; g.T = 0\} \subset \G$.  Using its unique element $g$ with $g(0) = 1$ we define the map
\[
\A \longrightarrow G 
\qquad
T\longmapsto g(1)
\]
which is surjective, has kernel $\Gn$ and thus gives an isomorphism $\A/\Gn \simeq G$.
We observe that the action of $\G$ on $\A$ descends to an action of $\G/\Gn\simeq G\times G$ on $G$.
\par
We now view the vector space $\A$ as an infinite-dimensional manifold with a natural flat metric. Using  an invariant scalar product $\langle\cdot,\cdot\rangle$  on $\g$, we have the norm  $\|t\|^2 = \int^1_0\langle t(s), t(s)\rangle ds$ for tangent vectors $t\in T_T\A$.
$\G$ and $\Gn$ act isometrically on $\A$, so we get a metric on $\A/\Gn$ as follows: at every point $m\in\A$ the tangent space of $\A$ splits into a subspace tangent to the orbit of the $\Gn$ action and its orthogonal complement, $T_m\A = T_m\A^\| \oplus T_m\A^\bot$. We define the metric at $[m] \in \A/\Gn$ via the identification of $T_{[m]}\A/\Gn$ with $T_m\A^\bot$.
\par
As $\A/\Gn \simeq G$, we now have a metric on $G$.
It is bi-invariant, since the action of $\G$ induces  an isometric action of $G\times G$.

\exercise{Compute the metric directly as the quotient metric.}

\subsection{Complexify}
We now complexify the construction of the previous section. This obviously should lead to a complexification of the Lie group $G$. 
\par
First we complexify the space of paths in the Lie algebra $\g$:
\[
\AC = \{T_0+iT_1 : [0,1] \longrightarrow \g\otimes\cx \}.
\]
$\AC$ is an infinite dimensional hermitian vector space with the hermitian metric
\begin{equation}
\|(t_0,t_1)\| = \int_0^1\|t_0(s)\|^2 + \|t_1(s)\|^2 ds.
\label{metric}\end{equation}
We write the corresponding symplectic form \[
\omega\left((t_0,t_1),(t_0',t_1')\right) = \int_0^1 \langle t_0,t_1'\rangle - \langle t_1,t_0'\rangle
\]
as $\omega = \int^1_0 dT_0 \wedge dT_1$. 

Extend the action of $\G$ on $\A$ to $\AC$ as follows:
\begin{equation}
\label{complexaction}
g.(T_0,T_1) = (gT_0g^{-1} - \dot{g}g^{-1}, gT_1g^{-1})
\end{equation}
This action again preserves all the structures, i.e. the metric and the symplectic form.

\claim{The action is \emph{Hamiltonian}, i.e. there is an equivariant  map $\mu:\AC\rightarrow (\Lie G)^*$ such that for any $\rho\in\Lie \Gn$ and any $v\in T\AC$
\[
\omega(\rho^*,v) = \langle d\mu(v), \rho\rangle
\qquad
\mathrm{where}\]
where $\rho^\ast$ is the vector field induced by $\rho$, i.e. $ 
\left. \rho^*\right|_m = \left. \frac{d}{d\epsilon}\left( \exp(\epsilon\rho).m\right) \right|_{\epsilon=0}$.
The map $\mu$ is called a \emph{moment map}.
}

To prove the claim, consider $\rho\in\Lie\Gn$, i.e. a path with $\rho(0) = 0 = \rho(1)$. Then, putting $g=\exp(\epsilon\rho)$ in \eqref{complexaction} and differentiating with respect to $\epsilon$, we obtain $\rho^* =\bigl([\rho, T_0] - \dot{\rho}, [\rho, T_1]\bigr)$. We now compute:
\begin{eqnarray*}
\omega\left(\rho^*, (t_0, t_1)\right) 
& = & \int_0^1 \left\langle [\rho, T_0] -\dot{\rho}, t_1 \right\rangle -\langle[\rho, T_1],t_0\rangle \\
\textrm{\small (integration by parts)} & = & - \langle \rho, t_1\rangle \Bigr|_0^1 + \int_0^1\langle \rho, \dot{t_1}\rangle + \langle [\rho, T_0], t_1\rangle - \langle [\rho, T_1], t_0\rangle \\
& = & \int_0^1\left\langle \rho, \dot{t_1} +[T_0, t_1]+ [t_0, T_1] \right\rangle\\
& = & \int_0^1 \Bigl\langle \rho, d(\underbrace{\dot{T_1} + [T_0, T_1]}_\mu) (t_0, t_1) \Bigr\rangle.
\end{eqnarray*}
Hence there is a moment map $\mu(T_0,T_1) = \dot{T_1}+[T_0, T_1]$, proving the claim.

Now recall the \emph{K\"ahler reduction}: let $(M,\langle\cdot,\cdot\rangle, I)$ be a K\"ahler manifold with its metric and complex structure and let $\mu:M\rightarrow \fH^*$ be a moment map for a free isometric and holomorphic action of a Lie group $H$ on $M$. If $c\in \fH$ is a fixed point of the co-adjoint action, then the quotient $\mu^{-1}(c)/H$ is again a K\"ahler manifold, called the \emph{K\"ahler quotient} and denoted by $M/\!/\!_c H$.
\par
We can apply this construction in our infinite-dimensional setting to $M=\AC$, $H=\Gn$, $c=0$. Since the action of $\Gn$ extends to a global action of $\Gn^\cx$, we can identify the K\"ahler quotient with a quotient of an open subset of $\AC$ by $\Gn^\cx$. This open subset is the union of \emph{stable} $\Gn^\cx$-orbits, i.e. those that meet $\mu^{-1}(0)$.  We shall see shortly that all $\Gn^\cx$-orbits are stable and so the same procedure as in the real case gives:
\begin{equation}
\AC/\!/\!_0 \,\Gn \simeq \AC/\Gn^\cx \simeq G^\cx.\label{complex}\end{equation}

What is the symplectic form on $G^\cx$? Let us analyse the K\"ahler quotient construction. The level set  $\mu^{-1}(0) $ of the moment map is given by the equation:
\begin{equation}
 \dot{T_1} = [T_1, T_0].
\label{level}\end{equation}
As in the previous section, quotienting by $\Gn$ is equivalent to finding a gauge transformation $g\in\G$ with $g(0) =1$ and sending $T_0$ to $0$. Equation \eqref{level} is gauge-invariant and hence $g$ makes $T_1$ constant. We obtain a map
\begin{equation}
\left( T_0(s), T_1(s)\right) \longmapsto \left( g(1), T_1(0)\right)
\label{sympl}\end{equation}
which gives an isomorphism 
\[
\mu^{-1}(0) / \Gn \simeq G \times \g \simeq T^*G,
\]
where the last isomorphism identifies $\g$ with right-invariant $1$-forms. $T^*G$, being a cotangent bundle, has a canonical symplectic form and an easy calculation shows that the symplectic form on $\mu^{-1}/\Gn$ coincides with the one on $T^\ast G$.
\par
We have now identified the K\"ahler manifold $\mu^{-1}(0) / \Gn$ with (possibly an open subset of) $G^\cx$ as a complex manifold, and with $T^\ast G$ as a symplectic manifold. To connect these two, we look again at the isomorphism \eqref{complex}, obtained by finding a complex gauge transformation $\tilde{g}$, $\tilde{g}(0) =1$, sending $T_0+iT_1$ to $0$. If $T_0+iT_1$ is already in $\mu^{-1}(0)$, i.e. it satisfies \eqref{level}, then we can find $\tilde{g}$ in two stages as follows. First, we find a real gauge transformation $g(s)$ sending $T_0$ to $0$ (with $g(0)=1$). This means that $g.(T_0+iT_1)=iT_1(0)$. This gives the map \eqref{sympl}. Now the transformation $\exp(isT_1(0))$ makes $g.(T_0+iT_1)$ identically zero and hence $\tilde{g}(s)=\exp(isT_1(0))g(s)$. Therefore $\tilde{g}(1)= \exp(iT_1(0))g(1)$ and so the K\"ahler metric on $G^\cx$ arises from the identification $G^\cx\simeq T^\ast G\simeq G\times \g$ given by the polar decomposition.
\par
This argument shows also that all $\Gn^\cx$-orbits are stable, as we have just identified the K\"ahler quotient with whole $G^\cx$ and not just with its open subset.
\par
We remark that the passage from $G$ with its bi-invariant metric to $T^\ast G\simeq G^\cx$ with its K\"ahler structure is an example of the \emph{adapted complex structure} construction (cf. \cite{LS}).

\subsection{Quaternionise}
\label{quaternionise}
Analogously to complexifying $\A$, we can quaternionise it:
\[
\AH = \{T_0 + i T_1 + j T_2 + k T_3 : [0,1] \longmapsto \g\otimes\quat \}
\]
On this space the natural $L^2$-metric (analogous to \eqref{metric}) is K\"ahler for three anti-commuting complex structures $I_1$, $I_2$, $I_3$ given by right multiplication by  $i$, $j$ and $k$. Thus we get three symplectic forms $\omega_1$, $\omega_2$, $\omega_3$, which can be explicitly computed. E.g. $\omega_1 = \int_0^1 dT_0 \wedge dT_1 - dT_2\wedge dT_3$\footnote[1]{Recall that $\int_0^1dT_0\wedge dT_1$ is a shorthand for $\omega\left((t_0,t_1),(t_0',t_1')\right) = \int_0^1 \langle t_0,t_1'\rangle - \langle t_1,t_0'\rangle$.}. In other words, $\AH$ is \emph{hyperk\"ahler}.

As before $\G$ acts on $\AH$ preserving the metric and the symplectic forms by extending the actions we had for the real and complexified cases:
\[
g.(T_0,T_1,T_2,T_3) = (gT_0g^{-1} - \dot{g}g^{-1}, gT_1g^{-1}, gT_2g^{-1}, gT_3g^{-1}). 
\]
The action of the subgroup $\Gn$ is free and Hamiltonian for all $3$ symplectic forms, giving us \emph{three} moment maps:
\begin{eqnarray}
\label{Nahm}
\mu_1 & = & \dot{T_1} + [T_0,T_1] - [T_2, T_3] \nonumber\\
\mu_2 & = & \dot{T_2} + [T_0,T_2] - [T_3, T_1] \\
\mu_3 & = & \dot{T_3} + [T_0,T_3] - [T_1, T_2]. \nonumber
\end{eqnarray}
We can perform the \emph{hyperk\"ahler reduction}, i.e consider the quotient
\[
\mu^{-1}_1(0) \cap \mu^{-1}_2(0) \cap \mu^{-1}_3(0) / \Gn
\]
which, as we shall show in a moment, is a hyperk\"ahler manifold.  The equations obtained by setting $\mu_i=0$ in \eqref{Nahm} are known as the 
\emph{Nahm equations} \cite{Nahm}.
\par
To see that we obtain a hyperk\"ahler manifold,   consider the complex-valued $2$-form $\omega_2+ i\omega_3$ which is a holomorphic $2$-form for the complex structure $I_1$. With respect to this complex symplectic form we get a holomorphic complex moment map 
\[
\mu_2+i\mu_3 : \AH \longrightarrow \Lie(\Gn)\otimes \cx\]
and so  
$(\mu_2+i\mu_3)^{-1} (0)$ is 
an $I_1$-complex submanifold of $\AH$, in particular K\"ahler. Using this we can identify the hyperk\"ahler quotient $
N = \mu^{-1}_1(0) \cap \mu^{-1}_2(0) \cap \mu^{-1}_3(0) / \Gn$ with the K\"ahler quotient of 
$(\mu_2+i\mu_3)^{-1} (0)$ by $ \Gn$. Doing this for the complex structures $I_2,I_3$ (and corresponding holomorphic $2$-forms) shows that $N$ is a hyperk\"ahler manifold.
\par
To identify $N$ let us go through this procedure in greater detail. We compute:
\begin{eqnarray*}
\mu_2+i\mu_3 &= & (\dot{T_2} + i \dot{T_3}) + [T_0, T_2+iT_3] -[T_3, T_1] -i[T_1,T_2] 
\\ & =& \frac{d}{ds}(T_2+iT_3)+ [T_0 - i T_1, T_2+iT_3].
\end{eqnarray*}
Just as in the real and complex cases, we now use the element $g\in \G^\cx$ with $g(0)=1$  which sends the first component $T_0 - iT_1$ to $0$ and $T_2 + iT_3$ to the constant $(T_2+iT_3)(0)$, giving a biholomorphism\footnote[1]{Again this requires proving that all $\Gn^\cx$-orbits on $(\mu_2+i\mu_3)^{-1} (0)$ are stable. An argument can be found in \cite{Don}: solving the remaining Nahm equation $\mu_1=0$ corresponds to finding a path in $G^\cx/G$ stationary under certain positive Lagrangian. Donaldson's argument also shows that any $\Gn^\cx$ orbit on $(\mu_2+i\mu_3)^{-1} (0)$ meets $\mu_1=0$ in a \emph{unique} $\Gn$-orbit - something that requires a proof in the infinite-dimensional setting.}
\[
(N, I_1) \longrightarrow G^\cx \times \overline{\g^\cx}
\qquad
(T_0,T_1,T_2,T_3) \longmapsto \left(g(1), (T_2 +i T_3)(0) \right),
\]
where $\overline{\g^\cx}$ denotes $\g^\cx$ with the opposite complex structure. Thus $(N, I_1)$ is biholomorphic to a cotangent bundle:
\[
(N, I_1) \simeq  
G^\cx \times \overline{ \g^\cx} \simeq 
T^{0,1}G \simeq
T^*G^\cx.
\]
The same construction can be done for the other complex structures $I_2$ and $I_3$.
\begin{remark} The quaternionisation procedure of this section does not yield a Lie group for the simple reason that $\g\otimes\quat$ is not a Lie algebra. Nevertheless there is an algebraic structure (Kronecker product) on $\g\otimes \quat$ obtained the same as for $\g\otimes \cx$: the Lie bracket in the first factor and the field (or skew-field) product in the second factor. The Nahm equations are closely related to this and  it has always been my feeling that there should exist an algebraic structure (non-associative product?) on $N\simeq T^\ast G^\cx$ reflecting the Kronecker product on $\g\otimes \quat$.
\end{remark}

\subsection{K\"ahler potentials}

First recall that on hyperk\"ahler manifolds we don't just have three complex structures at our disposal but a whole 2-sphere of them which is the unit sphere in the imaginary quaternions.
Our hyperk\"ahler manifold $N$ obtained in the previous section as the moduli space of solutions to Nahm's equations on the interval $[0,1]$ has a special property: all these complex structures are equivalent. In fact there is an isometric $SO(3)$-action on $N$ which induces the $SO(3)$-action on the $2$-sphere of complex structures. A matrix $A=[a_{ij}]\in SO(3)$ acts by
\[
(T_0,T_1,T_2,T_3) \longmapsto (T_0, A
\left(
\begin{array}{c} T_1 \\ T_2 \\ T_3 \end{array}
\right) )
= (T_0, \sum a_{1j}T_j, \sum a_{2j}T_j, \sum a_{3j}T_j).
\]

Knowing this $SO(3)$-action explicitly (as an action on $T^\ast G^\cx$, say) is of course
equivalent to knowing the hyperk\"ahler structure and it appears to be beyond reach. Instead, consider the $S^1 \subset SO(3)$ which fixes the complex structure $I_1$ (and $\omega_1$). It is given by
\[
(T_0,T_1,T_2+iT_3) \longmapsto
(T_0, T_1, e^{i\theta} (T_2 +iT_3))
\]
and the vector field induced by the action is $\theta^* = (0,0,-T_3 + iT_2)$. Compute the moment map for $\omega_1$:
\[
i_{\theta^*}\omega_1 
= \int_0^1 \langle T_3, dT_3\rangle - \left( - \langle dT_2, T_2\rangle \right) 
= \frac{1}{2} d \int_0^1 \|T_2\|^2 + \|T_3\|^2
\]
Hence the moment map for this action of $S^1$ is $\mu_{S^1} = ( \|T_2\|^2 + \|T_3\|^2)/2$. It turns out also to be a \emph{K\"ahler potential} for the complex structure $I_2$, i.e. $dI_2 d \mu_{S^1} = \omega_2$:
\begin{eqnarray*}
dI_2 d \mu_{S^1} 
&=& dI_2 \int_0^1 \langle T_3, dT_3\rangle + \langle T_2 , dT_2\rangle \\
&=& d\int_0^1 \langle T_3, -dT_1\rangle + \langle T_2, -dT_0\rangle \\
&=& \int_0^1 (dT_1\wedge dT_3 + dT_0 \wedge dT_2 ) = \omega_2
\end{eqnarray*}
This is indeed a general fact: an isometric circle action on a hyperk\"ahler manifold which preserves one of the complex structures and rotates the others will give a K\"ahler potential for another \cite{HKLR}.
\subsection{Spectral curves}
Following Hitchin, we can say something about the K\"ahler potential found in the previous section in the case of $G=SU(k)$.
\par

Consider again the level set for the complex moment map $\mu_2+i\mu_3$. It is given by the equation
\[
\dot{(T_2 +iT_3)} +  [ \underbrace{T_0 -iT_1}_\alpha,\underbrace{T_2 +iT_3}_\beta ] = 0
\qquad \textrm{i.e.} \qquad
 \dot{\beta} = [\beta, \alpha]
\]
which is a \emph{Lax equation} and thus implies that $\beta(t)$ lies in a fixed adjoint orbit of $G^\cx$.

If $G=SU(k)$ the spectrum of $\beta$ is constant and consists of $k$ points in $\cx$ but it depends on the complex structure we are using. Thus, for any complex structure, i.e. any point in $\proj^1$, we get $k$ points in $\cx$. This is equivalent to getting $k$ points in $T_\zeta\proj^1$ for every $\zeta\in\proj^1$, meaning that we are getting a $k$-fold ramified covering of $\proj^1$ in $T\proj^1$:  an algebraic curve $S$, called a \emph{spectral curve}. To understand $S$, notice that Nahm's equations are equivalent to 
$\dot{\beta}(\zeta) = [\beta(\zeta), \alpha(\zeta)]$ for all $\zeta\neq \infty$, where $\beta(\zeta) = \beta + (\alpha + \alpha^*)\zeta - \beta^*\zeta^2$ and $\alpha(\zeta) = \alpha -\beta^*\zeta$. Therefore $S$ is described as a compactification of the spectrum of $\beta(\zeta) $ for all $\zeta\neq\infty$.

For $G= SU(k)$, we have $S=\{(\eta,\zeta); \det(\eta \cdot 1 - \beta(\zeta)) = 0 \}$ where $\zeta\in\cx$ is an affine coordinate on $\proj^1$ and $\eta$ is a fibre coordinate of $T\proj^1$. Hence
\[
S = \{(\zeta,\eta);\eta^k + a_1(\zeta) \eta^{k+1} + \cdots + a_k(\zeta) = 0 \}
\]
where the $a_j$ are polynomials of degree $2k$.

Next try to get Nahm's equations back from the curve $S$. A theorem by Beauville \cite[1.4]{Beau1} gives that
\[
 \{A(\zeta)=A_1 + A_2 \zeta + A_3 \zeta^2 ; A_i \in \gl (k,\cx), \det(\eta\cdot 1 - A(\zeta)) = 0 \} / GL(k,\cx)
 \simeq
 J^{g-1}(S) -\Theta,
 \]
 where $g$ is the genus of $S$, $J^{g-1}(S)$ the Jacobian of line bundles of degree $g-1$ on $S$ and $\Theta$ the theta-divisor. The idea of the proof is to consider the sequence of sheaves on $T\proj^1$
\begin{equation}
0\longrightarrow
\pi^*\sO(-2)^k \longrightarrow
\pi^*\sO^k \longrightarrow
E\longrightarrow 
0,
\label{E}\end{equation}
where  the first map is $\eta\cdot 1 -A_1 -A_2\zeta - A_3\zeta^2$  and $E$ is the cokernel. $E$ is supported on $S$ and if $S$ is smooth, then $E$ is a line bundle of degree  $g+k-1$. After tensoring \eqref{E} with $\pi^\ast\sO(-1)$, we conclude that $E(-1)$ is a line bundle of degree $g-1$ with no global sections and, hence,
$E(-1)\in J^{g-1}(S) -\Theta$.
\par
Solutions to Nahm's equations for $U(k)$ correspond to a linear flow on $J^{g-1}(S)$ in a tangent direction corresponding to a fixed line bundle $L$ on $T\proj^1$ with $c_1(L)=0$. It has transition functions $e^{\eta/\zeta}$ and  induces a line bundle of degree $0$ on $S$ and thus gives us a vector field on $J^{g-1}(S)$.

\begin{theorem} (Hitchin, \cite[Prop 5.2]{Hit})  For a smooth $S$, let $\vartheta$ be the Riemann theta function on $J^{g-1}(S)\simeq J(S)$, where the isomorphism is given by $E\mapsto E(-k+2)$. Then the K\"ahler potential $\mu_{S^1}$ for $\omega_2$ on $T^*SL_n(\cx)$ is given by
\[
\mu_{S^1} = 
\frac{\vartheta'(a)}{a} - \frac{\vartheta'(b)}{b} + \frac{2\vartheta^{(N+2)}}{(N+1)(N+2)\vartheta^N(0)} - \frac{1}{6} a_2(1),
\]
where $a,b\in J^{g-1}(S)$ correspond to the triples $\bigl(T_1(0), T_2(0), T_3(0)\bigr)$ and  $\bigl(T_1(1), T_2(1)$,\linebreak$T_3(1)\bigr)$,  $'$ is the derivative in $L$ direction, $a_2(\zeta)$ is the second coefficient in $\eta^k + a_1(\zeta)\eta^{k-1} + a_2(\zeta)\eta^{k-2} + \cdots $ and $N = k(k^2-1)/6$.\end{theorem}

 With this we get a K\"ahler potential in terms of theta functions. It is non-algebraic, but this is only one problem. We now move to a setting where things are algebraic, and yet the hyperk\"ahler structure is hardly better understood.

\section{Orbits}

\subsection{Adjoint orbits and the `baby Nahm equation'}
Consider again the `baby Nahm equation', i.e. the Lax equation
\begin{equation}
\label{babynahm}
\dot{T_1} = [T_1,T_0]
\end{equation}
this time on a half-line:  $T_0, T_1: [0,+\infty) \rightarrow \g$. Assume that for $s\rightarrow\infty$, $T_0(s)$ converges  exponentially fast to $0$ and $T_1(s)$ converges exponentially fast to a fixed element $\tau\in\g$. For simplicity assume that $\tau$ is \emph{regular}, i.e. its centraliser $Z[\tau]$ is a Cartan subalgebra $\fH\subset\g$.

Let $\sN_\tau$ be the space of solutions to \eqref{babynahm} satisfying these boundary conditions and define
\[
\G = \{g:[0,\infty) \longrightarrow G : \lim_{t\rightarrow \infty}g(s) \in  \exp(\fH), \ \textrm{converging exponentially fast}\}
\]
Also define the subgroup $\Gn = \{g\in\G : g(0) =1\}$. While these groups and solutions have different boundary conditions from those we used in the previous sections, we can define the same action as in (\ref{complexaction}):
\[
T_0 \longmapsto gT_0 g^{-1} -\dot{g}g^{-1}
\qquad\textrm{and}\qquad
T_1 \longmapsto gT_1g^{-1}
\]
For every $(T_0,T_1)\in\sN_\tau$ we find  an element $g\in\Gn$ such that
\[
g.(T_0, T_1) = (0,T_1(0)) 
\qquad \textrm{where} \qquad
T_1(0) = g(\infty) \tau g(\infty)^{-1}.
\]
Thus we get that the quotient by the action is the \emph{adjoint orbit} of $\tau$: $\sN_\tau/\Gn \simeq \sO_\tau$. As in the previous sections this is a K\"ahler quotient.
Its symplectic form is the Kostant-Kirillov-Souriau form on $\sO_\tau$: identifying $T_x\sO_\tau $ with $ \{[\rho,x] : \rho\in\g \}$ we have:
\[
\omega([\rho,x], [\rho',x]) = \langle[\rho,\rho'], x\rangle.
\]
 To find out the complex structure of our K\"ahler quotient, let, again, $\alpha = T_0 - iT_1$ and rewrite the `baby Nahm equation' as $\dot{\alpha+\alpha^*} = [\alpha^*,\alpha]$, where $\alpha$ converges exponentially fast to $i\tau$ as $t\rightarrow\infty$. Hence we can find $g\in\G^\cx$ such that $\alpha = g(i\tau)g^{-1} - \dot{g}g^{-1}$.

Note that $g$ is not unique: Given $g$ satisfying the equation, we can replace it by $g'(t) = g(t) e^{i\tau t} p e^{-i\tau t} $, where $p$ is  constant, to give us the same $\alpha$. However, the resulting $g'(t)$ will  converge to an element of the Cartan  subalgebra $\fH$ if and only if $\mathrm{ad}(i\tau)$ acts on $p$ with non-negative eigenvalues. That is, if $p$ is in the Borel subgroup $B$ determined by $\mathrm{ad}(i\tau)$.

With this, the map $\alpha\mapsto g(0)$ gives a biholomorphism between $\sN_\tau/\Gn$ and $G^\cx/B$: the complex structure of $\sO_\tau$ is that of a generalised flag variety.

\exercise{Adapt the construction to non-regular $\tau$.}

\subsection{Adjoint orbits and Nahm's equations}
We now quaternionise to return to Nahm's equations
\begin{eqnarray*}
\dot{T_1} + [T_0,T_1] - [T_2, T_3] &=&0\nonumber\\
\dot{T_2} + [T_0,T_2] - [T_3, T_1] &=&0\\
\dot{T_3} + [T_0,T_3] - [T_1, T_2] &=&0. \nonumber
\end{eqnarray*}
While the equations are the same as in section \ref{quaternionise},  this time our solutions are defined on the half-line $[0,+\infty)$. Assume that we have exponentially fast convergence
\[
T_0(s) \longrightarrow 0
\qquad \textrm{and} \qquad
T_i(s) \longrightarrow \tau_i
\qquad \textrm{for} \qquad s\longrightarrow \infty,
\]
where the $\tau_i$ are fixed. The Nahm equations imply that $[\tau_i,\tau_j]=0$ for $i,j=1,2,3$. Thus if $\tau_i$ are regular for $i=1,2,3$, then they lie in a common Cartan subalgebra $\fH$.  Actually we will assume a bit less than regularity for each $\tau_i$, namely that the centraliser of the triple $Z(\tau_1, \tau_2,\tau_3)=Z(\tau_1)\cap Z(\tau_2)\cap Z(\tau_3)$ is a Cartan subalgebra $\fH$.

Let $\sN_{\ttt}$ denote the set of solutions of Nahm's equations with the given boundary conditions. Then the quotient $\sM_{\ttt} = \sN_{\ttt}/\Gn$ is a hyperk\"ahler manifold \cite{Kron}.
\par
To find out what are the complex structures of $\sM_{\ttt}$ we write, as before, $\alpha = T_0 - iT_1$ and $\beta = T_2 +iT_3$. We can write\footnote[1]{Once again, this isomorphism depends on the fact, proved in \cite{Kron}, that all $\Gn$-orbits are stable.}
\[
(\sM_{\ttt},I_1) = \{(\alpha, \beta) : \dot{\beta} = [\beta, \alpha]\}/\Gn^\cx.
\]
If $\tau_2+i\tau_3$ is regular in $\Lie{G}^\cx$, then the map $(\alpha, \beta) \longmapsto \beta(0)$ gives a biholomorphism
\[
(\sM_{\ttt},I_1) \longrightarrow \sO_{\tau_2+i\tau_3}
\]
where $\sO_{\tau_2 + i \tau_3}$ denotes the adjoint $G^\cx$-orbit of $\tau_2+i\tau_3$. The idea for proving this is to find a $g\in\G^\cx$ such that
\[
(\alpha,\beta) = \left(g(i\tau_1) g^{-1} - \dot{g}g^{-1}, g(\tau_2+i \tau_3) g^{-1}\right).
\]
We conclude that the generic complex structure of $\sM_{\ttt}$ is that of a complex adjoint orbit. The other complex structures are those of some holomorphic bundles over generalised flag manifolds:

\example{If $\tau_2=\tau_3 = 0$, then $\tau_1$ is regular and we can write (cf. end of the previous subsection)
\[
(\alpha,\beta) = \left(g(i\tau_1) g^{-1} - \dot{g}g^{-1}, g(\exp(i\tau_1 t)n \exp(i\tau_1 t)) g^{-1}\right),
\]
where $n$ is in the nilradical $\fN$ of $\fB$, the Borel subalgebra given by $\ad (i\tau_1)$. Thus we have a map
\[
\left(\sM_{\tau_1,0,0}, I_1\right)\longrightarrow G^\cx\times_B \fN \simeq T^\ast\bigl(G^\cx/B\bigr)
\qquad
(\alpha,\beta) \longmapsto g(0)
\]
which is biholomorphic.
}

What about other orbits? Notice that there are other non-constant `standard' solutions to Nahm's equations. For example, consider $\su(2)$ with the basis $\{e_1, e_2,e_3\}$ and relations $[e_1, e_2] = -e_3, ...$.  If $\sigma:\su(2)\rightarrow\g$ is a homomorphism then
\begin{equation}
\label{nil}
T_i(s) = \sigma(e_i)/(s+1) \qquad \textrm{for } i = 1,2,3
\end{equation}
is a solution to Nahm's equations. The Jacobson-Morozow theorem gives us a one-to-one correspondence between nilpotent orbits in $\g^\cx$ and homomorphisms $\su(2) \rightarrow \Lie{G}$ up to conjugation. If we consider solutions to Nahm's equations asymptotic to \eqref{nil}, we also get a hyperk\"ahler metric, this time on an appropriate nilpotent orbit \cite{Kron2}.
 
Most generally, we fix $\tau_1,\tau_2,\tau_3 \in \fH$ and $\sigma: \su(2) \rightarrow Z(\ttt)$ to get the set
\[
\sM_{\ttt;\sigma} = \{ \textrm{solutions asymptotic to } \tau_i + \frac{\sigma(e_i)}{s+1} \textrm{ for } i=1,2,3\}/\Gn,
\]
where both ``asymptotic" and $\Gn$ need to be defined carefully \cite{Biq,Kov}.
\par
This is again a hyperk\"ahler manifold.
Its generic complex structure is that of the orbit of
\[
(\underbrace{\tau_2+i\tau_3}_{\textrm{\footnotesize semisimple}}) + (\underbrace{\sigma(e_2) +i \sigma(e_3)}_{\textrm{\footnotesize nilpotent}})
\]
if $Z(\tau_2,\tau_3) = Z(\tau_1, \tau_2,\tau_3)$. The complex-symplectic form $\omega_2+i\omega_3$ is the Kostant-Kirillov-Souriau form of the complex co-adjoint orbit.

With respect to other complex structures $\sM_{\ttt;\sigma}$ is a bundle over a generalised flag manifold.
\par
Even this is not the most general construction of hyperk\"ahler metric on (co)-adjoint orbits -- see \cite{Santa-Cruz}. Nevertheless the manifolds $\sM_{\ttt;\sigma}$ are probably the most general \emph{algebraic} hyperk\"ahler structures on coadjoint orbits. They are also most general \emph{geodesically complete} hyperk\"ahler structures on coadjoint orbits:
\begin{theorem}[\cite{Biel}] If $M$ is a complete hyperk\"ahler manifold with a hyperk\"ahler action of  a compact semisimple Lie group $G$ such that for one complex structure $G^\cx$ acts locally transitively, then $M$ is equivariantly isomorphic to one of the $\sM_{\ttt}$.
\end{theorem}

The idea for proving the theorem is to take the moment maps $\mu_1, \mu_2, \mu_3$ for the $G$-action on $M$ and consider the gradient flow maps $m(t)$ of $\|\mu_1\|^2$, where $I_1$ is a complex structure for which $G^\cx$ is locally transitive. Then the $T_i$ given by  $T_i(t) = \mu_i(m(t))$ satisfy Nahm's equations.

\subsection{Explicit description of the hyperk\"ahler metrics on adjoint orbits?}

Biquard \cite{Biq2} shows that the hyperk\"ahler structure of $\sM_{\ttt;\sigma}$ is algebraic. Nevertheless, there is no algebraic construction of these structures, except for a few special cases (\cite{BG,KS1,KS2}). One would like, at the very least, to relate the hyperk\"ahler structure of $\sM_{\tau_1,0,0}$ to the K\"ahler structure of $\sO_{\tau_1}$. This has been done by Biquard and Gauduchon \cite{BG}, but only when $\sO_{\tau_1}$ is a Hermitian symmetric space.
\par
The other possibility (for $G=SU(k)$) is suggested by Hitchin \cite{Hit}, analogous to his approach to $T^\ast SL(k,\cx)$ described in section 1.5. For $\sM_{\ttt}$, we have a fixed spectral curve $S$ given by the equation
\begin{equation}\det\bigl(\eta\cdot 1 -(\tau_2+i\tau_3) -2i\tau_1\zeta +(-\tau_2+i\tau_3)\zeta^2\bigr)=0.\label{sing}\end{equation}
$S$ is reducible, in fact a union of rational curves. When $\tau_2=\tau_3=0$, we have, as in section 1.4, a K\"ahler potential and once again we can try to describe it in terms of theta functions. This time, however, the theta functions are rational \cite{Mum}.
\par
While nobody succeeded yet in describing the hyperk\"ahler structure of $\sM_{\tau_1,0,0}$ using this approach, thinking along these lines leads to an amusing observation. Recall that if the curve \eqref{sing} has only nodes, then $J^{g-1}(S)$ has a canonical compactification $\overline{J^{g-1}(S)}$ obtained by adding invertible sheaves of degree $g^\prime-1$ and semistable multidegree on partial normalisations $S^\prime$ of $S$. The definition of a semistable multidegree is one which allows us to extend the notion of the theta-divisor to
$\overline{J^{g-1}(S)}$: semistable multi-degrees are those for which the usual definition of the theta-divisor actually gives a divisor \cite{Beau1}.
\par
Observe also that the curve $S$ is invariant under the antiholomorphic involution
$$\zeta\mapsto -1/\bar{\zeta}, \quad \eta\mapsto -\eta/\bar{\zeta}^2.$$
This induces an antiholomorphic involution on $H^1(S,\sO_S)$ and we can speak of \emph{real} line bundles of degree zero. Similarly we define a line bundle $L$ of degree $g-1$ to be real if $L(-k+2)$ is real (here $G=SU(k)$). This extends to the compactified Jacobian and for any $U\subset \overline{J^{g-1}(S)}$ we write $U_\oR$ for the corresponding real sheaves.
\par
We have:
\begin{proposition} Let $\tau_2+i\tau_3$ be a regular semisimple element of $SL(k,\cx)$ and let $\sO_{\tau_2+i\tau_3}$ denote its adjoint orbit. Choose a  $\tau_1$ in the same Cartan subalgebra as $\tau_2$ and $\tau_3$  so that the curve $S$ defined by \eqref{sing} has only nodes as singularities. Then there is a canonical algebraic isomorphism between $\sO_{\tau_2+i\tau_3}/SU(k)$ and the closure in $\left(\overline{J^{g-1}(S)}-\Theta\right)_{\oR}$ of a connected component of $\left(J^{g-1}(S)-\Theta\right)_{\oR}$.
\end{proposition}

\section{Applications to symmetric pairs and real orbits}

\subsection{Kostant-Sekiguchi correspondence}
Let $(G,K)$ be a compact symmetric pair, i.e. there exists an orthogonal decomposition
\[
\g = \kk\oplus\fM \quad \textrm{ with} \quad
[\fK, \fM] \subset \fM, \enskip [\fM,\fM] \subset \fK.
\]
Then $\g^\ast = \fK + i \fM$ is also a Lie algebra and if $G^\ast$ is the corresponding Lie group, then $(G^\ast,K)$ is the dual symmetric pair.
\example{ $G=SU(n)$ and $K=S(U(p)\times U(q))$ with $p+q=n$. $\fM$ consists of the two off-diagonal $p\times q$- and $q\times p$-blocks in the algebra of skew-hermitian matrices and $SU(n)^\ast = SU(p,q)$, i.e the group preserving an indefinite hermitian form.}
 \example{ $G=SU(n)$ and $K=SO(n)$.
 Then $\fM$ is the space of imaginary symmetric matrices. It follows that $SU(n)^\ast = SL(n,\oR)$.}
\medskip

The Kostant-Sekiguchi correspondence tells us that there is a one-to-one correspondence between nilpotent orbits of $G^\ast$ and $K^\cx$-orbits of nilpotent elements in $\fM^\cx \subset \Lie{G}^\cx$. This somewhat mysterious fact has been given a natural interpretation in terms of Nahm's equations by Vergne \cite{Vergne}. We shall explain this now.

\subsection{ ($\g,\kk$)-valued solutions to Nahm's equations}
Let $\sM$ be any  moduli space of $\g$ valued solutions to Nahm's equations given as a quotient of the space $\sN$ of solutions  with fixed boundary conditions by an appropriate gauge group $\Gn$. Following Saksida \cite{Saks}, we consider the subset $\sN^{\gk}$ of \emph{$\gk$-valued solutions}, i.e. those solutions in $\sN$ with $T_0(s), T_1(s) \in \kk$ and $T_2(s), T_3(s) \in \fM$ for all $s$. We define their moduli space as
\[
\sM^{\gk} = \sN^{\gk} /\{g\in \Gn; g(s) \in K\enskip \textrm{for all $s$}\}
\]
The usual $\g$-valued solutions to Nahm's equations can be thought of as $\sM^{(\g\oplus\g, \g)}$ .

Recall that having chosen a complex structure we can write the moduli space of solutions as the quotient of the level set of the complex moment map $\dot{\beta} = [\beta, \alpha]$ by the complexified Lie group.
Look again at the $\gk$-valued solutions and choose complex structures $I_1, I_3$. This gives
\[
\beta_1 = T_2 + i T_3, \qquad \alpha_1 = T_0 - iT_1,
\qquad
\beta_3 = T_1+ i T_2, \qquad \alpha_3 = T_0 -i T_3,
\]
satisfying $\dot{\beta_i} = [\beta_i,\alpha_i]$ for $i=1,3$, and where, for all $s$, $\alpha_1(s)\in \fK^\cx, \beta_1(s)\in\fM^\cx, \alpha_3(s)\in  \g^\ast, \beta_3(s) \in  \g^\ast$. Let us also define two subgroups of the gauge group $\Gn^\cx$:
$$ \sK^\cx_0=\{g\in \Gn^\cx; g(s) \in K^\cx\enskip \textrm{for all $s$}\},$$
$$ \sG^\ast_0=\{g\in \Gn^\cx; g(s) \in G^\ast\enskip \textrm{for all $s$}\}.$$

Then we  have identifications
\begin{eqnarray*}
\sM^{\gk}
&=& \sN^{\gk} / \sK_0 \\
&\stackrel{\phi_1}{\simeq} &
\{\dot{\beta_1} = [\beta_1, \alpha_1] : \alpha_1\in\kk^\cx, \beta_1 \in \fM^\cx\} / \sK_0^\cx\\
&\stackrel{\phi_3}{\simeq}&
\{\dot{\beta_3} = [\beta_3, \alpha_3] :  \alpha_3\in\g^\ast, \beta_3 \in \g^\ast\} / \G^\ast_0.
\end{eqnarray*}
The proof of these identifications relies on  solving the third Nahm equation: this, as for the usual Nahm's equations (see the footnote on p. 5), is equivalent to finding a stationary path for a positive geodesically convex potential  in negatively curved spaces $K^\cx\!/K$ and $G^\ast\!/K$.
\par
We now consider the spaces $\sM^{\gk}$ for the moduli spaces of Nahm's equations considered in these lectures: i.e. $T^\ast G^\cx$ and adjoint orbits.
\par
For the moduli space $N$ defining a hyperk\"ahler structure on $T^\ast G^\cx$, $N^{\gk}$ is isomorphic via $\phi_1$ to $K^\cx \times \fM^\cx$ which is a complex submanifold of $G^\cx \times\overline{\g^\cx}$ and via $\phi_3$ to $G^\ast\times \g^\ast \simeq T^\ast G^\ast$. With this we get a canonical $K$-invariant K\"ahler metric on $T^\ast G^\ast$.

\subsection{Real orbits} It is more interesting to consider $\sM^{\gk}_\ttts$, where $\sM_\ttts$ are the moduli spaces described in section 2, defining hyperk\"ahler structures on adjoint orbits of $G^\cx$.
Assume that $I_1, I_3$ are generic so that
\[
(\sM_\ttts, I_1) \simeq \sO_{v_1}
\qquad \textrm{with} \qquad
v_1 = \tau_2 +i \tau_3 + \sigma(e_2) + i \sigma(e_3)
\]\[
(\sM_\ttts, I_3) \simeq \sO_{v_3}
\qquad \textrm{with} \qquad
v_2 = \tau_1 +i \tau_2 + \sigma(e_1) + i \sigma(e_2).
\]
We assume that the restricted moduli space $\sM_\ttts^{\gk}$ is not empty, i.e. $\tau_1\in \fK$ and $\tau_2,\tau_3\in \fM$.  Then, restricting the isomorphism $(\alpha,\beta) \mapsto \beta(0)$ gives, as in \cite{Kron, Biq, Kov}:
\[
\phi_1: \sM_\ttts^{\gk} \stackrel{\sim}{\longrightarrow} \sO_{v_1} \cap \fM^\cx
\qquad \textrm{and} \qquad
\phi_3: \sM_\ttts^{\gk} \stackrel{\sim}{\longrightarrow} \sO_{v_3} \cap \g^\ast.
\]
Thus we have a diffeomorphism
\begin{equation} \sO_{v_1}\cap\fM^\cx \stackrel{\sim}{\longrightarrow} \sO_{v_3}\cap \g^\ast. \end{equation}
 In the case $\tau_1=\tau_2=\tau_3 = 0$ this is the Vergne isomorphism \cite{Vergne}.
 \par
 Observe that
$\sO_{v_3} \cap \g^\ast$ is a finite union of $G^\ast$-orbits and $\sO_{v_1}\cap \fM^\cx$ is $K^\cx$-invariant. To show that $\sO_{v_1}\cap \fM^\cx$ is a union of $K^\cx$-orbits we need

\begin{lemma}(cf. \cite{Bryl}) $K^\cx$ acts transitively on connected components of $\sO_{v_1}\cap\fM^\cx$.\end{lemma}
\begin{proof}
We need to show that the action of $K^\cx$ on $\sO_{v_1}\cap\fM^\cx$ is infinitesimally transitive. Let $x\in\sO_{v_1}\cap\fM^\cx$ and consider
\[
T_x\sO_{v_1} = \{[\rho,x] : \rho \in \g^\cx \}
\simeq [\kk^\cx, x] \oplus [\fM^\cx, x]
\subset \fM^\cx \oplus \fK^\cx.
\]
With this we have that $T_x\sO_{v_1} \cap \fM^\cx \simeq [\fK^\cx, x] $ and so $K^\cx$ acts locally transitively as required.\end{proof}

We can extend the Vergne diffeomorphism and the Kostant-Sekiguchi correspondence to non-nilpotent orbits as follows:
\begin{proposition} Let $(G,K)$, $\g=\kk\oplus \fM$, be a compact symmetric pair. Let $\tau_1\in \fK$, $\tau_2,\tau_3\in\fM$ be such that their $\g$-centralisers  satisfy: $Z(\tau_1)\supset Z(\tau_2)\cap Z(\tau_3)$, $Z(\tau_3)\supset Z(\tau_1)\cap Z(\tau_2)$. Then there is a one-to-one correspondence between $G^\ast$-orbits whose closure contains the orbit of $\tau_1+i\tau_2$ and $K^\cx$-orbits in $\fM^\cx$ whose closure contains the $K^\cx$-orbit of $\tau_2+i\tau_3$. Moreover, for each pair of corresponding orbits, there exists a canonical $K$-equivariant diffeomorphism between them.
\end{proposition}

In particular, setting $\tau_1=\tau_3=0$, we obtain:
\begin{corollary} Let $(G,K)$, $\g=\kk\oplus \fM$, be a compact symmetric pair. There is a $1-1$ correspondence between adjoint $G^\ast$-orbits of elements whose semisimple part lies in $i\fM$ and $K^\cx$-orbits of elements of $\fM^\cx$ whose semisimple part lies in $\fM$. Moreover, for each pair of corresponding orbits, there exists a canonical $K$-equivariant diffeomorphism between them.
\end{corollary}

We finish with a simple example of the Vergne diffeomorphism. Let $\sO$ be the non-zero nilpotent orbit in $\fS\fL_2(\cx)$. $\sO$ can be identified with $\bigl(\cx^2-
 0 \bigr)/ \Z_2$  using the map
\begin{equation}(u,v) \longmapsto \left(\begin{array}{cc}uv & u^2 \\ -v^2 & -uv \end{array}\right).\label{nilp}
\end{equation}
We have the usual identification of  $\cx^2$ with $\quat=\{x_0+ix_1+jx_2+kx_3\}$ (for the complex structure $i$) by letting $u= x_0+ix_1$ and $ v = x_2+ix_3$. The hyperk\"ahler structure on $\sO$ is then that of $\bigl(\quat- 0 \bigr)/ \Z_2$.
With respect to the complex structure corresponding to $j$, the map \eqref{nilp} is
\begin{equation}(u,v) \longmapsto \left(\begin{array}{cc}(u-i\bar{v})(v+i\bar{u}) & (u-i\bar{v})^2 \\ -(v+i\bar{u})^2 & (i\bar{v}-u)(v+i\bar{u}) \end{array}\right).\label{nilp2}
\end{equation}

Consider now $\sO \cap \fS\fL_2(\oR) = \{u^2,v^2,uv \in \oR\}$. It consists of two orbits $\sO_+$ and $\sO_-$ of $\fS\fL_2(\oR)$ given respectively by $u^2+v^2 > 0$ and  $ u^2+v^2 < 0$.  The condition $u^2+v^2>0$ is equivalent to $x_1=0$ and $x_3=0$, while $u^2+v^2<0$ is equivalent to $x_0=0$ and $x_2=0$. A direct computation show that the map \eqref{nilp2} sends $\sO_+$ and $\sO_-$ to the sets of matrices of the form
$$\left(\begin{array}{cc} ib & b \\ b & -ib \end{array}\right)\quad\textrm{for $\sO_+$},\qquad
\left(\begin{array}{cc} -ib & b \\ b & +ib \end{array}\right)\quad\textrm{for $\sO_-$}.$$
These are precisely the two non-zero nilpotent orbits of $K^\cx=SO(2,\cx)$ on $\fM^\cx=\{$symmetric matrices$\}$.


\end{document}